\documentclass[a4paper,11pt]{article} 
\usepackage{float}
\usepackage{tikz-qtree}
\usetikzlibrary{arrows.meta}
\usepackage[a4paper, margin=.9in]{geometry}
\usepackage{amsmath,amsthm,amsfonts,amssymb,array}
\usepackage[square,sort,comma,numbers]{natbib}
\usepackage{subeqnarray}
\usepackage{blkarray}
\usepackage[utf8]{inputenc}
\usepackage[english]{babel}
\usepackage{mathrsfs}
\usepackage{colortbl}
\usepackage{graphicx,epsfig}
\newcounter{myeqno}
\usepackage{color}
\usepackage{tcolorbox}
\usepackage{ upgreek }
\usepackage{ mathrsfs }
\usepackage{relsize}
\usepackage{multirow}
\definecolor{shadecolor}{gray}{0.75}

\theoremstyle{definition}

\newtheorem{proposition}{Proposition}[section]

\newtheorem{theorem}{Theorem}[section]
\newtheorem{corollary}{Corollary}[section]
\newtheorem{problem}{Problem}[section]
\newtheorem{lemma}{Lemma}[section]

\usepackage[misc,geometry]{ifsym}
\date{}
\usepackage[colorlinks=true]{hyperref}
\hypersetup{urlcolor=blue, citecolor=blue}
\newlength{\defbaselineskip}
\newcommand{\setlinespacing}[1]%
{\setlength{\baselineskip}{#1 \defbaselineskip}}

\begin{document}

		\title{\textbf{Distance Seidel matrix of a connected graph}}

\author{Haritha T$^1$\footnote{harithathottungal96@gmail.com},  Chithra A. V$^1$\footnote{chithra@nitc.ac.in}
	\\ \\ \small 
 1 Department of Mathematics, National Institute of Technology Calicut,\\\small Calicut-673 601, Kerala, India\\ \small}

\maketitle	

	\thispagestyle{empty}
		\begin{abstract}
    For a connected graph $G$, we present the concept of a new graph matrix related to its distance and  Seidel matrix, called distance Seidel matrix $\mathcal{D}^S(G)$. Suppose that the eigenvalues of $\mathcal{D}^S(G)$ be $\partial_{1}^{S}(G) \geq \cdots \geq \partial_{n}^{S}(G).$ In this article, we establish a relationship between distance Seidel eigenvalues of a graph with its distance and adjacency eigenvalues. We characterize all the connected graphs with $\partial_{1}^{S}(G)= 3.$ Also, we determine different bounds for the distance Seidel spectral radius and distance Seidel energy. We study the distance Seidel energy change of the complete bipartite graph due to the deletion of an edge. Moreover,   we obtain the distance Seidel spectra of different graph operations such as join, cartesian product, lexicographic product, and unary operations like the double graph and extended double cover graph. We give various families of distance Seidel cospectral and distance Seidel integral graphs as an application.
\end{abstract}

{Keywords:Distance matrix, Eigenvalues, Seidel matrix, Energy, Spectral radius, Cospectral graphs.}
\section{Introduction}
   Let $G$ be an undirected simple connected graph with vertex set $V(G)= \{v_r: 1\leq r \leq n\}$ and an edge set $E(G)=\{e_t: 1\leq t\leq m\}$, then \textit{adjacency matrix} $A(G)= (a_{rt})$ of $G$ is an $n\times n$ symmetric matrix, where $a_{rt}= 1$ if $v_r$ and $v_t$ are adjacent in $G$, zero otherwise. Let the eigenvalues of $A(G)$ ($A$-eigenvalues) be $\lambda_1\geq \lambda_2\geq \cdots \geq \lambda_n$, respectively. The complement $\Bar{G}$ of $G$ is
defined by taking $V(\Bar{G})= V(G)$ and making two vertices $u$ and $v$ adjacent in $\Bar{G}$ if and only if they are non-adjacent in $G.$ The adjacency matrix of $\Bar{G}$ is denoted by $\Bar{A}.$ In $1966$, Van Lint and Seidel defined the Seidel matrix of graph $G$ as $$\mathcal{S}(G)= J-I- 2A(G)= (s_{rt})_{n\times n}= \begin{cases}-1, &\text{if $v_r$ is adjacent to $v_t$,}\\0, &\text{if r=t}\\1, &\text{if $v_r$ is non adjacent to $v_t$.}\end{cases}$$ Then the eigenvalues of $\mathcal{S}(G)$ are said to be the Seidel eigenvalues of $G$ ($\mathcal{S}$-eigenvalues), denoted by $\lambda^{S}_1\geq \lambda^{S}_2\geq\cdots \geq \lambda^{S}_n.$ The collection of all $\mathcal{S}$-eigenvalues of $G$ is said to be the Seidel spectrum of $G.$ In \cite{haemers2012seidel}, Haemers defined the Seidel energy $E_{\mathcal{S}}(G)$ of $G$ as the sum of absolute values of all $\mathcal{S}$-eigenvalues of $G$.  For more results on the energy and spectrum of Seidel matrix of graphs, see \cite{oboudi2016energy,akbari2020some}. A graph is said to be \textit{integral} if all its eigenvalues are integers.  For a partitioned matrix, we denote its \textit{equitable quotient matrix} by $M$ \cite{brouwer2011spectra}. Also all entry one matrix is denoted by $J.$ Let the degree of vertex $v_r$ in $G$ be $d_r.$ The \textit{complete graph}, \textit{path}, \textit{wheel graph}, \textit{cycle}, \textit{complete multipartite graph} on $q$ parts ($2\leq q \leq n-1$), and \textit{star graph} are denoted by $K_n$, $P_n$, $W_n$, $C_n$, $K_{n_1, \ldots, n_q}$, and $S_n$ respectively, and when $q=2, \; K_{n_1, n_2}$ is said to be the complete bipartite graph. The \textit{complete split graph} $CS(n, p)$ with parameters $n, p \;(p\leq n)$, is a graph of order $n$ containing a clique of order $p$ and an independent set of order $n-p$ in which every vertex of the independent set is adjacent to each vertex of the clique. The \textit{friendship graph} $f_n$ is the join of $K_1$ with $nK_2$.\\
\par The distance $d_G(v_r,v_t)$ or $d_{rt}$ between any two vertices $v_r$ and $v_t$ in $V(G)$, is the length of the shortest path connecting them in $G$ \cite{buckley1990distance}, also $d_G(v_r, v_r)=0$ for any $v_r\in V(G)$. The \textit{diameter} $\tilde{\omega}$ of $G$ is $\tilde{\omega}= max\{d_G(v_r,v_t): v_r,v_t \in V(G)\}$. The \textit{distance matrix} $\mathcal{D}= \mathcal{D}(G)$ of $G$, is a square symmetric matrix defined by
$\mathcal{D}(G)= (d_{rt})_{n\times n}= d_G(v_r, v_t).$ Let $\partial_1\geq \partial_2\geq \cdots \geq \partial_s$ denote the distinct eigenvalues of $\mathcal{D}(G)$, called the distance eigenvalues of $G$ ($\mathcal{D}$-eigenvalues). Then the collection of all $\mathcal{D}$-eigenvalues of $G$ is called the distance spectrum $Spec_{\mathcal{D}}(G)$ ($\mathcal{D}$-spectrum) of $G.$ The studies on spectral properties of $\mathcal{D}(G)$ can be seen in \cite{aouchiche2014distance}. Graphs having all its $\mathcal{D}$-eigenvalues as integers are called \textit{distance integral graphs}. The \textit{distance energy} $E_{\mathcal{D}}(G)$ \cite{indulal2010distance} of $G$ is defined as $E_{\mathcal{D}}(G)= \sum_{r=1}^{n}\partial_r.$  \\
\par The \textit{Wiener index} $\mathcal{W}(G)$ \cite{wiener1947structural} of a graph $G$ is a topological index defined as $$\mathcal{W}(G)= \frac{1}{2}\sum_{{\substack{v_r, v_t\in V(G)}}} d_{G}(v_r, v_t).$$ The \textit{transmission} $Tr(v_r)$ of a vertex $v_r$  is the sum of all distances from $v_r$ to all other vertices in $G.$ A graph $G$ satisfying $Tr(v_r)= k$ for every vertices $v_r$ in $G$ is called a $k$-transmission regular graph. In \cite{aouchiche2013two}, the authors defined the Laplacian and signless Laplacian for a distance matrix called the \textit{distance Laplacian matrix} and the \textit{distance signless Laplacian matrix}.\\

\par Motivated by the above works we introduce distance Seidel matrix, which is similar to the Seidel matrix in graph theory, and it is defined as
$$\mathcal{D}^S(G)= J-I-2\mathcal{D}= \begin{cases}-1, &\text{if $v_r$ is adjacent to $v_t$,}\\1-2d_{G}(v_r, v_t), &\text{otherwise.}\end{cases}$$
  \\
Let the eigenvalues of $\mathcal{D}^S(G)$ be $\partial_1^S\geq \partial_2^S \geq \cdots \geq \partial_n^S$, called the distance Seidel eigenvalues ($\mathcal{D}^S$-eigenvalues), and the set of all $\mathcal{D}^S$-eigenvalues of $G$ is said to be the distance Seidel spectrum $Spec_{\mathcal{D}^S}(G)\; (\mathcal{D}^S$-spectrum). The distance Seidel spectral radius is defined as,
$\rho^S(G)= max_{i}(|\partial_i^{S}|).$ The distance Seidel energy $E_{\mathcal{D}^S}(G)$ of a graph $G$ is defined as $$E_{\mathcal{D}^S}(G)= \sum_{r=1}^{n}|\partial_r^S|.$$

\par This paper is organized as follows. In Section $2$, we give some known results. In Section $3$, we establish the relation of $\mathcal{D}^S$-eigenvalues with $\mathcal{D}$-eigenvalues and $A$-eigenvalues and study some other spectral properties of $\mathcal{D}^S$-eigenvalues. Also we characterize connected graphs with $\mathcal{D}^S$-eigenvalue $3$. Section $4$, deals with different bounds of distance Seidel spectral radius. In Section $5$, we give the $\mathcal{D}^S$-energy of various graphs and obtain bounds for the $\mathcal{D}^S$-energy. Moreover we provide, one open problem based on the $\mathcal{D}^S$-energy of complete bipartite graph by the deletion of an edge. In Section $6$, we determine $\mathcal{D}^S$-spectrum of graphs formed by using different graph operations. Finally, in Section $7$, we obtain different families of distance Seidel cospectral and distance Seidel integral graphs.
\section{Preliminaries}
We give some known results in this section, which are used in the main results.
\begin{lemma} \cite{davis1979circulant}
Let $B= \begin{bmatrix}B_0& B_1\\ B_1& B_0\end{bmatrix}_{2\times 2}$ be a block symmetric matrix. Then the eigenvalues of $B_0+B_1$ together with those of $B_0-B_1$ forms the eigenvalues of $B.$
\end{lemma}
\begin{lemma}{(Cauchy Interlace Theorem)} \cite{bapat2010graphs}
Let $P$ be a Hermitian matrix of order $n$ with eigenvalues $\tilde{\lambda}_{1}(P) \geq \tilde{\lambda}_{2}(P) \geq \cdots \geq \tilde{\lambda}_{n}(P)$, and let $Q$ be a principal submatrix of $P$ of order $m$ with eigenvalues $\tilde{\lambda}_{1}(Q) \geq \tilde{\lambda}_{2}(Q) \geq \cdots \geq \tilde{\lambda}_{n}(Q).$ Then $\tilde{\lambda}_{n-m+r}(P) \leq \tilde{\lambda}_{r}(Q) \leq \tilde{\lambda}_{r}(P)$ for $1 \leq r \leq m.$
\end{lemma}
\begin{proposition} \cite{so1994commutativity}
Let $P$ and $Q$ be any two $n\times n$ Hermitian matrices and $R= P+Q.$ Then
\begin{itemize}
    \item [(i)]$\tilde{\lambda}_r(R)\leq \tilde{\lambda}_t(P)+ \tilde{\lambda}_{r-t+1}(Q),\; 1\leq t \leq r \leq n$,
    \item[(ii)]$\tilde{\lambda}_r(R)\geq \tilde{\lambda}_t(P)+\tilde{\lambda}_{r-t+n}(Q), \; 1\leq r \leq t \leq n$, where $\tilde{\lambda}_r(P)$ is the $r^{th}$ largest eigenvalue of $P$.
\end{itemize}
\end{proposition}
Let $B$ be a real $n\times n$ matrix with non negative entries. If $B^k> 0$ for some $k$, then $B$ is said to be primitive and it is \textit{irreducible} if there exist a $k$ such that $(B^k)_{rt}> 0$ for all $r, t.$ 
\begin{theorem} \cite{indulal2008distance}
Let $G$ be a $k$-transmission regular graph of order $n$ with $\mathcal{D}$-eigenvalues $k= \partial_1\geq \partial_2\geq \cdots \geq \partial_n$. Then there exist a polynomial $p(x)$ such that $p(\mathcal{D})= J,$
$$p(x)= \frac{n(x-\partial_2)\ldots(x-\partial_n)}{(k-\partial_2)\ldots(k-\partial_n)}.$$
\end{theorem}
\begin{lemma} \cite{zhou2007largest}
Let $C= (c_{rt})_{n\times n}$ be an irreducible symmetric matrix of order $n\geq 2$ with row sums $C_1, C_2, \ldots, C_n$ and largest eigenvalue $\rho(C)$, then\\
$$\sqrt{\frac{\sum_{r=1}^{n}C_{r}^2}{n}}\leq \rho(C)\leq max_{1\leq r\leq n}\sum_{t=1}^{n}c_{rt}\sqrt{\frac{C_t}{C_r}}$$
with equality if and only if $C_1= C_2= \dots= C_n$ or there is a permutation matrix $Q$ such that $Q^TCQ= \begin{bmatrix}0& B\\ B^T& 0\end{bmatrix}.$
\end{lemma}
\begin{theorem} \cite{lin2013distance}
Let $\mathcal{D}$ be the distance matrix of a connected graph $G$, then $\partial_n(\mathcal{D})= -2$ with multiplicity $n-q$ if and only if $G$ is a complete $q$-partite graph for $2\leq q \leq n-1.$ 
\end{theorem}

 \section{Distance Seidel eigenvalues of a graph and its properties}
This section deals with the basic properties of distance Seidel eigenvalues of $G$ and their relation between adjacency and distance eigenvalues.
\begin{proposition}
If $G$ is a $k$-transmission regular graph of order $n$, then its $\mathcal{D}^S$-spectrum is given by $$Spec_{\mathcal{D}^S}(G)= \begin{pmatrix}n-1-2k& -1-2\partial_r\\ 1& 1\end{pmatrix},\; \text{where}\; r= 2, 3, \ldots, n.$$
\end{proposition}
\begin{proof}
The proof follows from  Theorem $2.1$
\end{proof}
\subsection{$\mathcal{D}^S$-eigenvalues of some standard graphs}

1. The complete graph
\\

$$Spec_{\mathcal{D}^S}(K_n)= \begin{pmatrix}1-n& 1\\ 1& n-1\end{pmatrix}.$$\\
\\
2. Cycle
\\
\par Using the $\mathcal{D}$-spectrum of $C_n$, we get the $\mathcal{D}^S$-spectrum of $C_n$ as follows:\\
For odd $n= 2k+1$, the $\mathcal{D}^S$-eigenvalues are $n-\frac{n^2}{2}-\frac{1}{2}$, $\left(-1+\frac{1}{2}\operatorname{sec}^2(\frac{\pi r}{n})\right)^2,\; r= 1,\ldots, k.$
For even $n= 2k$, the $\mathcal{D}^S$-eigenvalues are $n-1-\frac{n^2}{2}$, $(-1)^{k-1}$, $\left(-1+2\operatorname{csc}^2(\frac{\pi(2r-1)}{n})\right)^2,\; r= 1, \ldots, \left\lfloor\frac{k}{2}\right\rfloor$, and $1$ if $k$ is odd.\\
\\
3. The complete bipartite graph
\\
$$Spec_{\mathcal{D}^S}(K_{a,b})= \begin{pmatrix}3& \frac{3(2-a-b)\pm \sqrt{9a^2+9b^2-14ab}}{2}\\a+b-2& 1\end{pmatrix}.$$\\
\\
4. The star graph
\\
 $$Spec_{\mathcal{D}^S}(S_n)= \begin{pmatrix}3& \frac{6-3n\pm \sqrt{9n^2-32n+32}}{2}\\n-2& 1\end{pmatrix}.$$\\
\\
5. The complete graph due to an edge deletion
\\
$$Spec_{\mathcal{D}^S}(K_n-e)= \begin{pmatrix}3& 1& \frac{-n\pm \sqrt{n^2-4n+20}}{2}\\1& n-3&1 \end{pmatrix}.$$
\par Next we give the Wiener index $\mathcal{W}(G)$ of $G$ in terms of distance and distance Seidel eigenvalues of $G$.
\begin{proposition}
For a connected graph $G$ of order $n$,
$$\mathcal{W}(G)= \frac{1}{8}\left(n(n-1)-\sum_{r=1}^{n}{(\partial_{r}^{S}})^{2}+4\sum_{r=1}^{n}\partial_r^2\right).$$
\end{proposition}
\begin{proof}
We know that,
\begin{equation*}
    \begin{aligned}
    \sum_{r=1}^{n}{(\partial_{r}^{S}})^{2}&= trace((\mathcal{D}^S)^{2})\\
    &= \sum_{s=1}^{n}\sum_{t=1}^{n}(1-2d_{st})^2\\
    &= 2\sum_{s<t}^{n}(1-2d_{st})^2\\
    &= n(n-1)-8\sum_{s<t}^{n}d_{st}+8\sum_{s<t}^{n}d_{st}^2\\
    &= n(n-1)-8\mathcal{W}(G)+4\sum_{r=1}^{n}\partial_r^2.
    \end{aligned}
\end{equation*}
  
Therefore, $\mathcal{W}(G)= \frac{1}{8}\left(n(n-1)-\sum_{r=1}^{n}{(\partial_{r}^{S}})^{2}+4\sum_{r=1}^{n}\partial_r^2\right).$
\end{proof}

\begin{theorem}
Let $G$ be a connected graph of order $n$ with $\tilde{\omega}\leq 2$. If $\lambda_1, \lambda_2, \ldots, \lambda_n$ and $\partial_1^S, \partial_2^S, \ldots, \partial_n^S$ are the $A$-eigenvalues and $\mathcal{D}^S$-eigenvalues of $G$ respectively, then
$$3-3n+2\lambda_r\leq \partial_r^S\leq 3+2\lambda_r,\; 1\leq r \leq n.$$
\end{theorem}
\begin{proof}
By definition, $\mathcal{D}^S(G)= J-I-2\mathcal{D}= 3I-3J+2A$.\\
From  Proposition $2.1$, we have
 \begin{equation}\tag{1}
 \tilde{\lambda}_r(\mathcal{D}(G)^S)\leq \tilde{\lambda}_t(3I-3J)+\tilde{\lambda}_{r-t+1}(2A(G)),
    \end{equation}
 \begin{equation}\tag{2}
     \tilde{\lambda}_r(\mathcal{D}^S(G))\geq -\tilde{\lambda}_t(3J-3I)+\tilde{\lambda}_{r-t+n}(2A(G)).
 \end{equation}\\
Put $t=1$ in (1), we get
\begin{equation}\tag{3}
   \partial_r^S\leq 3+2\lambda_r,\;  r= 1, \ldots, n.
\end{equation}\\
Put $t=n$ in (2), we get
\begin{equation}\tag{4}
    \partial_r^S\geq 3-3n+2\lambda_r,\; r= 1, \ldots, n.
\end{equation}\\
From equations (3) and (4), we get the desired result.

\end{proof}
\begin{theorem}
If $\partial_1, \partial_2, \ldots, \partial_n$ and $\partial_1^S, \partial_2^S, \ldots, \partial_n^S$ are the $\mathcal{D}$-eigenvalues and $\mathcal{D}^S$-eigenvalues of $G$ respectively, then
$$-1-2\partial_{n-t+1}\leq \partial_t^S\leq n-1-2\partial_{n-t+1}, \; 1\leq t \leq n.$$
\end{theorem}
\begin{proof}
By applying   Proposition $2.1$, we get\\
 \begin{equation}\tag{5}
   \tilde{\lambda}_r(J-I)\leq \tilde{\lambda}_t(\mathcal{D}^S(G))+\tilde{\lambda}_{r-t+1}(2\mathcal{D}(G)),
 \end{equation}\\
 \begin{equation}\tag{6}
   \tilde{\lambda}_r(J-I)\geq \tilde{\lambda}_t(\mathcal{D}^S(G))+\tilde{\lambda}_{r-t+n}(2\mathcal{D}(G)).
  \end{equation}
Now let $r=n$ in (5), we get\\
\begin{equation}\tag{7}
    \partial_t^S\geq -1-2\partial_{n-t+1}.
\end{equation}\\
Put $r=1$ in (6), we get\\
\begin{equation}\tag{8}
    \partial_t^S\leq n-1-2\partial_{n-t+1}.
\end{equation}\\
By combining (7) and (8), we get the desired result.
\end{proof}
\begin{theorem}
For a connected graph $G$ of order $n$, $\partial_{1}^{S}= 1$ if and only if $G\cong K_n.$
\end{theorem}
\begin{proof}
If $G\cong K_n$, then trivially, $\partial_{1}^{S}= 1.$ Conversely, let $\partial_{1}^{S}(G)= 1$. If $\tilde{\omega}\geq 2$, then $\mathcal{D}^S(P_3)$ is a principal submatrix of  $\mathcal{D}^S(G)$, with $\partial_{1}^{S}(P_3)= 3.$ By Lemma $2.2$, we get $\partial_{n-2}^{S}(G) \leq \partial_{1}^{S}(P_3) \leq \partial_{1}^{S}(G)$, which is a contradiction (since $\partial_{1}^{S}(P_3)= 3$). That is,  $G\cong K_n.$
\end{proof}
\begin{lemma}
Let $G= K_{n_1, \ldots, n_q}$ $(2\leq q \leq n-1)$ of order $n$ with $n_1\geq n_2\geq \cdots \geq n_q.$ Then $\partial_{1}^{S}(G)= 3$ of multiplicity $n-q.$
\end{lemma}
\begin{proof}
Since $G= K_{n_1, \ldots, n_q},\; \tilde{\omega}= 2$, then $\mathcal{D}^S(G)= I-J-2\Bar{A}$. By  Proposition $2.1$, $\tilde{\lambda}_{1}(\mathcal{D}^S(G))\leq \tilde{\lambda}_{1}(I-J)+\tilde{\lambda}_{1}(-2\Bar{A})= 3.$ Next we show that there exist an eigenvalue which equals $3$. Then the distance Seidel matrix of $G$ is\\

\[\mathcal{D}^S(G)=\left(\begin{array}{ccccc}
3 (I_{n_{1}, n_{1}}-J_{n_{1}}) & J_{n_{1}, n_{2}} & J_{n_{1}, n_{3}} & \cdots & J_{n_{1}, n_{q}} \\
J_{n_{2}, n_{1}} & 3( I_{n_{2}, n_{2}}- J_{n_{2}}) & J_{n_{2}, n_{3}} & \cdots & J_{n_{2}, n_{q}} \\
\vdots & \vdots & \ddots & \cdots & \vdots \\
J_{n_{q-1}, n_{1}} & J_{n_{q-1}, n_{2}} & J_{n_{q-1}, n_{3}} & \cdots & J_{n_{q-1}, n_{q}} \\
J_{n_{q}, n_{1}} & J_{n_{q}, n_{2}} & J_{n_{q}, n_{3}} & \cdots & 3 (I_{n_{q}, n_{q}}- J_{n_{q}})
\end{array}\right).\]
\\
$\begin{aligned}
&\operatorname{det}\left(xI-\mathcal{D}^S\right)\\
&= (x-3)^{n-q}\left(\begin{array}{ccccc}
x+3n_{1}-3 & -n_{2} & -n_{3}& \cdots & -n_{q} \\
-n_{1} & x+3n_{2}-3 & -n_{3}& \cdots & -n_{q} \\
\vdots & \vdots & \ddots & \cdots & \vdots \\
-n_{1} & -n_{2} & -n_{3}& \cdots & -n_{q} \\
-n_{1} & -n_{2} & -n_{3}& \cdots & x+3n_{q}-3
\end{array}\right) .
\end{aligned}$
\\
Thus the multiplicity of the $\mathcal{D}^S$-eigenvalue $3$ is at least $n-q\geq 1.$ Clearly $\mathcal{D}^S(G)= I-J-2\Bar{A}$.Then applying  Proposition $2.1$, we get $\tilde{\lambda}_{n}(I-J)+ \tilde{\lambda}_{r}(-2\Bar{A}) \leq \tilde{\lambda}_{r}(\mathcal{D}^S) \leq \tilde{\lambda}_{1}(I-J)+\tilde{\lambda}_{r}(-2\Bar{A}), \; 2\leq r \leq n.$ Here $\Bar{A}$ is the adjacency matrix of $\Bar{G}$ which is the union of complete graphs $K_{n_1}, \ldots , K_{n_q}$ so that $\lambda_{r}(\Bar{A})= n_{r}-1$ for $1\leq r \leq q.$ Therefore, $3-n-2n_{r} \leq \partial_{r}^{S} \leq 3-2n_{r}$ for $r= 2, \ldots, q.$ Thus $\partial_{r}^{S}\leq 1$, since $n_{r}\geq 1$ for  $2\leq r \leq q$ and $\partial_{n}^{S}\leq 3-n.$ Therefore the multiplicity of $3$ is $n-q.$
\end{proof}
\begin{theorem}
Let $G= K_{n_{1}, \ldots , n_{q}}\; (2\leq q \leq n-1)$ of order $n$ with $n_1\geq n_2\geq \cdots \geq n_q.$ Then the characteristic polynomial of $\mathcal{D}^S(G)$ is,
$$p(x)= (x-3)^{n-q}\left(\prod_{r=1}^{q} (x-3+4n_{r})-\sum_{r=1}^{q}n_{r}\prod_{t=1, t \neq r}^{q}(x-3+4n_{t})\right).$$
\end{theorem}
\begin{proof}
From the proof of  Lemma $3.1$, we have
\begin{equation*}
    \begin{aligned}
    &\operatorname{det}\left(xI-\mathcal{D}^S\right)\\
    &= (x-3)^{n-q}\left(\begin{array}{ccccc}
x-3+3n_{1} & -n_{2} & -n_{3}& \cdots & -n_{q} \\
-n_{1} & x-3+3n_{2} & -n_{3}& \cdots & -n_{q} \\
\vdots & \vdots & \ddots & \cdots & \vdots \\
-n_{1} & -n_{2} & -n_{3}& \cdots & -n_{q} \\
-n_{1} & -n_{2} & -n_{3}& \cdots & x-3+3n_{q}
\end{array}\right)
\end{aligned}
\end{equation*}
By applying elementary row operations, we get
\begin{equation*}
    \begin{aligned}\operatorname{det}\left(xI-\mathcal{D}^S\right)
&= (x-3)^{n-q}\left(\prod_{r=1}^{q} (x-3+4n_{r})-\sum_{r=1}^{q}n_{r}\prod_{t=1, t \neq r}^{q}(x-3+4n_{t})\right).
    \end{aligned}
\end{equation*}

\end{proof}
Next we consider some forbidden subgraphs of $G$ which are used in the subsequent theorem. We call $F$ a forbidden subgraph of $G$ if $\mathcal{D}^S(F)$ is a proper principal submatrix of $\mathcal{D}^S(G)$ and $\partial_{1}^{S}(F)>3.$\\
 
 \hspace{1.5cm} 
 \begin{tikzpicture}[scale=0.4,inner sep=1pt]
\draw (0,0) node(1) [circle,draw,fill] {}
      (2,0) node(2) [circle,draw,fill] {}
      (1,2) node(3) [circle,draw,fill] {}
      (4,0) node(4) [circle,draw,fill] {}
      (8,0) node(5) [circle,draw,fill] {}
      (10,0) node(6) [circle,draw,fill] {}
      (12,0) node(7) [circle,draw,fill] {}
      (8,2) node(8) [circle,draw,fill] {}
      (10,2) node(9) [circle,draw,fill] {}
      (16,0) node(10) [circle,draw,fill] {}
      (18,0) node(11) [circle,draw,fill] {}
      (20,0) node(12) [circle,draw,fill] {}
      (22,0) node(13) [circle,draw,fill] {}
      (26,0) node(14) [circle,draw,fill] {}
      (28,0) node(15) [circle,draw,fill] {}
      (26,2) node(16) [circle,draw,fill] {}
      (28,2) node(17) [circle,draw,fill] {}
      (27,3) node(18) [circle,draw,fill] {};

\draw [-] (1) to (2) to (3) to (1);  
\draw [-] (2) to (4);
\draw[-] (5) to (6) to (7);
\draw[-] (5) to (8) to (9) to (6);
\draw[-] (10) to (11) to (12) to (13);
\draw[-] (14) to (15) to (17) to (18) to (16) to (14); 
\end{tikzpicture}

\noindent
Figure $3.1$. \hspace{0.8cm}$F_1$  \hspace{2.4cm} $F_2$ \hspace{3.2cm} $F_3$ \hspace{2.5cm} $F_4.$ \\
\\

 We have $\partial_{1}^{S}(F_1)= 3.78, \partial_{1}^{S}(F_2)= 5.97, \partial_{1}^{S}(F_3)= 5.82, \text{and}\; \partial_{1}^{S}(F_4)= 4.23.$ That is $\partial_{1}^{S}(F_{i})>3$, so $F_i$ is a forbidden subgraph of $G$ for $i= 1, 2, 3, 4$ and if $G$ contains $F_i$ as an induced
subgraph, then $\partial_{1}^{S}(G)>3.$
\begin{theorem}
Let $G$ be a connected graph of order $n$. Then $\partial_{1}^{S}=3$ with multiplicity $n-q$ if and only if $G$ is a complete $q$-partite graph where $q= 2, \ldots, n-1.$
\end{theorem}
\begin{proof}
If $G$ is a complete $q$-partite graph, then by Lemma $3.1$, $\partial_{1}^{S}= 3$ with multiplicity $n-q.$\\
Conversely, let $G$ be a connected graph with $\partial_{1}^{S}=3$ of multiplicity $n-q.$ \\
\textbf{Case 1.} $G$ is a tree.
\par If $\tilde{\omega} \geq 3$, then $\mathcal{D}^{S}(F_{3})$ is a proper principal submatrix of $\mathcal{D}^{S}(G)$, which is contradiction (since $\partial_{1}^{S}(F_3)>3$). Therefore, $G$ is of diameter $2$ with eigenvalue $3$, and so $G\cong K_{1,n-1}.$\\
\textbf{Case 2.} $G$ is not a tree.\\
By applying similar arguments as in the proof of Theorem $2.2$, we get the desired result.
\end{proof}
\section{Distance Seidel spectral radius and its bounds}
The distance Seidel spectral radius of $G$ is the largest eigenvalue of the irreducible matrix $-\mathcal{D}^S(G)$. If $e= uv$ is an edge of $G$ such that $G-e$ is connected, then $d_G(v_r, v_t)\leq d_{G-e}(v_r, v_t)$ for all $v_r, v_t\in V(G).$ Thus by Perron-Frobenius theorem $\rho^S(G)< \rho^S(G-e).$ From which we can say that the distance Seidel spectral radius is maximum for trees. Similarly, $\rho^S(G+e)< \rho^S(G),$ which implies the minimum distance Seidel spectral radius is attained by $K_n.$\\
This section gives different bounds for the distance Seidel spectral radius of a connected graph.
\begin{theorem}
Let $G$ be a connected graph of order $n\geq 2$ and for $1\leq r\leq n, \;\mathcal{D}^S_r= \sum_{{\substack{t= 1 \\ r \neq t}}}^{n}(2d_{rt}-1)$. Then\\
$$\rho^S(G) \leq max_{\substack{1\leq r\leq n}} \sum_{{\substack{t= 1 \\ r \neq t}}}^{n}(2d_{rt}-1)\sqrt{\frac{\mathcal{D}^S_t}{\mathcal{D}^S_r}}.$$
Equality holds if and only if $\mathcal{D}^S_1= \mathcal{D}^S_2= \dots= \mathcal{D}^S_n.$
\end{theorem}
\begin{proof}
Clearly $-\mathcal{D}^S(G)$ is an irreducible matrix having largest eigenvalue $\rho^S(G).$ Then by right inequality of Lemma $2.3$, we have\\
$$\rho^S(G) \leq max_{\substack{1\leq r\leq n}}\sum_{{\substack{t= 1 \\ r \neq t}}}^{n}(2d_{rt}-1)\sqrt{\frac{\mathcal{D}^S_t}{\mathcal{D}^S_r}}.$$
It follows that, equality holds if and only if $\mathcal{D}^S_1= \mathcal{D}^S_2= \dots= \mathcal{D}^S_n.$
\end{proof}
\begin{theorem}
Let $G$ be a connected graph of order $n\geq 2$ and $\mathcal{D}^S_r= \sum_{{\substack{t= 1 \\ r \neq t}}}^{n}(2d_{rt}-1)$ for $1\leq r \leq n.$ Then\\
$$\rho^S(G)\geq \sqrt{\frac{\sum_{r=1}^{n}(\mathcal{D}^S_r)^2}{n}}.$$
Equality holds if and only if $\mathcal{D}^S_1= \mathcal{D}^S_2= \dots= \mathcal{D}^S_n.$
\end{theorem}
\begin{proof}
From left inequality of Lemma $2.3$, we have $\rho^S(G)\geq \sqrt{\frac{\sum_{r=1}^{n}(\mathcal{D}^S_r)^2}{n}}$, and equality holds if and only if $\mathcal{D}^S_1= \mathcal{D}^S_2= \dots= \mathcal{D}^S_n.$
\end{proof}
\begin{theorem}
For a connected graph $G$ of order $n$, let $\Updelta$ and $\Updelta'$ denote the maximum degree and second maximum degree, respectively. Then\\
$$\rho^S(G)\geq \sqrt{(3n-2\Updelta-3)(3n-2\Updelta'-3)}.$$
Equality holds if and only if the graph $G$ is regular and $\tilde{\omega}\leq 2.$
\end{theorem}
\begin{proof}
Consider $-\mathcal{D}^S(G)= -(J-I-2D)$, which is an irreducible matrix whose largest eigenvalue is $\rho^S(G)$. Let $z=(z_1, z_2, \ldots, z_n)^T$ be the Perron eigenvector of $\mathcal{D}^S(G)$ corresponding to the eigenvalue $\rho^S(G)$, such that\\
$$z_{r}=\min _{t \in V(G)} z_{t} \quad \text { and } \quad z_{s}=\min _{\substack{t \in V(G) \\ t \neq r}} z_{t}.$$

Then from the equation $-\mathcal{D}^S(G)z= \rho^S(G)z$, we have
\begin{equation*}
    \begin{aligned}
    \rho^S(G)z_r&= \sum_{\substack{t= 1 \\ t \neq r}}^{n}(2d_{rt}-1)z_t\\
    &\geq d_rz_s+3(n-d_r-1)z_s= (3n-2d_r-3)z_s
    \end{aligned}
\end{equation*}

Similarly for the component $z_s$ we have\\
\begin{equation*}
    \begin{aligned}
    \rho^S(G)z_s&= \sum_{\substack{t= 1 \\ t \neq s}}^{n}(2d_{st}-1)z_t\\
    &\geq d_sz_r+3(n-d_s-1)z_r= (3n-2d_s-3)z_r.
    \end{aligned}
\end{equation*}

Combining these two inequalities, it follows\\
$$\rho^S(G)\geq \sqrt{(3n-2d_r-3)(3n-2d_s-3)} \geq \sqrt{(3n-2\Updelta-3)(3n-2\Updelta'-3)}.$$

First let  $\rho^S(G)= \sqrt{(3n-2\Updelta-3)(3n-2\Updelta'-3)}$ then we can see that $\tilde{\omega}\leq 2$ and all coordinates of $z$ are equal. Also we get $\rho^S(G)= 3n-2d_r-3= 3n-2d_s-3$, that is $d_r= d_s$. Therefore the graph $G$ is regular. Conversely assume that $G$ is a regular graph with $\tilde{\omega}\leq 2$. Then for $\tilde{\omega}=1$, we get the graph $K_n.$ If $\tilde{\omega}= 2$, we get $\rho^S(G)= 3n-2\Updelta-3.$

\end{proof}
\begin{theorem}
Let $G$ be a connected graph of order $n$ and diameter $\tilde{\omega}$, with minimum degree $\updelta$ and second minimum degree $\updelta'$, respectively. Then\\
{\scriptsize$\rho^S(G)\leq \sqrt{(-2\updelta(\tilde{\omega}-1)+\tilde{\omega}(1-\tilde{\omega})+2n\tilde{\omega}-n-1)(-2\updelta'(\tilde{\omega}-1)+\tilde{\omega}(1-\tilde{\omega})+2n\tilde{\omega}-n-1)}.$}\\
Equality holds if and only if the graph $G$ is regular with $\tilde{\omega}\leq 2.$
\end{theorem}
\begin{proof}
Let the Perron eigenvector corresponding to the largest eigenvalue $\rho^S (G)$ of $-\mathcal{D}^S(G)$ be $z= (z_1, z_2, \ldots, z_n)^T$ such that\\
$$z_{r}=\max _{t \in V(G)} z_{t} \quad \text { and } \quad z_{s}=\max _{\substack{t \in V(G) \\ t \neq r}} z_{t}.$$
Then from the equation $-\mathcal{D}^S(G)z= \rho^S(G)z$, we have
\begin{equation*}
    \begin{aligned}
    \rho^S(G) z_r&= \sum_{\substack{t= 1 \\ t \neq r}}^{n} (2d_{rt}-1)z_t\\
    &\leq d_rz_s+3z_s+\ldots+(2\tilde{\omega}-3)z_s+(2\tilde{\omega}-1)(n-1-d_r-(\tilde{\omega}-2))z_s\\
    &= (\tilde{\omega}-1)^2z_s-z_s+d_rz_s+(2\tilde{\omega}-1)(n-d_r-\tilde{\omega}+1)\\
    &= (-2d_r(\tilde{\omega}-1)+\tilde{\omega}(1-\tilde{\omega})+2n\tilde{\omega}-n-1)z_s.
    \end{aligned}
\end{equation*}

Similarly for the component $z_s$ we have\\
\begin{equation*}
    \begin{aligned}
    \rho^S(G) z_s&= \sum_{\substack{t= 1 \\ t \neq s}}^{n} (2d_{st}-1)z_t\\
    & \leq d_sz_r+3z_r+\ldots+(2\tilde{\omega}-3)z_r+(2\tilde{\omega}-1)(n-1-d_s-(\tilde{\omega}-2))z_r\\
    &= (-2d_s(\tilde{\omega}-1)+\tilde{\omega}(1-\tilde{\omega})+2n\tilde{\omega}-n-1)z_r.
    \end{aligned}
\end{equation*}
 
Combining these two inequalities, it follows\\
{\scriptsize\begin{equation*}
    \begin{aligned}
    \rho^S(G)&\leq \sqrt{(-2d_r(\tilde{\omega}  -1)+\tilde{\omega}(1-\tilde{\omega})+2n\tilde{\omega}-n-1)(-2d_s(\tilde{\omega}-1)+\tilde{\omega}(1-\tilde{\omega})+2n\tilde{\omega}-n-1)}\\
    \\
    &\leq \sqrt{(-2\updelta(\tilde{\omega}-1)+\tilde{\omega}(1-\tilde{\omega})+2n\tilde{\omega}-n-1)(-2\updelta'(\tilde{\omega}-1)+\tilde{\omega}(1-\tilde{\omega})+2n\tilde{\omega}-n-1)}.
    \end{aligned}
\end{equation*}}
\\
First let {\scriptsize$\rho^S(G)= \sqrt{(-2\updelta(\tilde{\omega}-1)+\tilde{\omega}(1-\tilde{\omega})+2n\tilde{\omega}-n-1)(-2\updelta'(\tilde{\omega}-1)+\tilde{\omega}(1-\tilde{\omega})+2n\tilde{\omega}-n-1)}$}, then we can see that all coordinates of $z$ are equal. Also, we get $\rho^S(G)= -2d_r(\tilde{\omega}-1)+\tilde{\omega}(1-\tilde{\omega})+2n\tilde{\omega}-n-1= -2d_s(\tilde{\omega}-1)+\tilde{\omega}(1-\tilde{\omega})+2n\tilde{\omega}-n-1$, that is $d_r= d_s$. For the case when $\tilde{\omega}\geq 3$, that is, for every vertex $r$, there is exactly one vertex $s$ which is at a distance $2$ from $r$, then $\tilde{\omega}$ of $G$ must be less than $4$. When $\tilde{\omega}= 3$ and equality holds, we can find a center vertex $s$ with $e(s)= 2$, then\\ 
\begin{equation*}
    \begin{aligned}
    \rho^S(G)z_s&= d_sz_s+(n-1-d_s)3z_s\\
    &= -2\updelta(3-1)+3(1-3)+2n\times3-n-1
    \end{aligned}
\end{equation*}\\
which implies $d_s= n-2$. That is, $G\cong P_4$, and the coordinates of Perron's eigenvector of  $D^S(P_4)$ are not all equal. Therefore, graph $G$ is regular with $\tilde{\omega}\leq 2$. Conversely let the graph $G$ is regular with $\tilde{\omega}\leq 2$. Then for $\tilde{\omega}=1$, we get the complete graph $K_n.$ For $\tilde{\omega}=2$, we get $\rho^S(G)= 3n-2\updelta-3.$
\end{proof}
From Theorems $4.3$ and $4.4$, we get the following corollary.
\begin{corollary}
For a connected graph $G$ of order $n$, denote the maximum degree, second maximum degree, minimum degree, and second minimum degree by $\Updelta$, $\Updelta'$, $\updelta$, and $\updelta'$ respectively. Then
$$\tilde{\alpha}(G)\leq \rho^S(G)\leq \tilde{\beta}(G).$$
where,
{\scriptsize\begin{equation*}
    \begin{aligned}
    &\tilde{\alpha}(G)= \sqrt{(3n-2\Updelta-3)(3n-2\Updelta'-3)},\\
    \\
    &\tilde{\beta}(G)= \sqrt{(-2\updelta(\tilde{\omega}-1)+\tilde{\omega}(1-\tilde{\omega})+2n\tilde{\omega}-n-1)(-2\updelta'(\tilde{\omega}-1)+\tilde{\omega}(1-\tilde{\omega})+2n\tilde{\omega}-n-1)}.
    \end{aligned}
\end{equation*}}

\end{corollary}

Now we present a lower bound for $\rho^S(G)$ for a bipartite graph $G$.
\begin{theorem}
Let $G$ be a bipartite connected graph with bipartition $V(G)= P\cup Q,\;|P|=p,\;|Q|=q,\;p+q= n.$ Let $\Updelta_P$ and $\Updelta_Q$ be the maximum degrees among vertices from $P$ and $Q$, respectively. Then\\
$$\rho^S(G)\geq \frac{3n-6+\sqrt{9n^2-36pq+4(5p-4\Updelta_Q)(5q-4\Updelta_P)}}{2}.$$
Equality holds if and only if $G= K_{p,q}$, or the graph $G$ is semi-regular with eccentricity $e(v_i)= 3$ for every $v_i\in V(G).$
\end{theorem}
\begin{proof}
Let $P= \{1, 2, \ldots, p\}$, $Q= \{p+1, p+2, \ldots, p+q\}$ and suppose that the Perron eigenvector of $-\mathcal{D}^S(G)$ corresponding to the eigenvalue $\rho^S(G)$ be $z= (z_1, z_2, \ldots, z_n)^T$ such that,\\
$$z_{r}=\min _{t \in P} z_{t} \quad \text { and } \quad z_{s}=\min _{t \in Q} z_{t}.$$
Then from the equation $-\mathcal{D}^S(G)z= \rho^S(G)z$, we have\\
\begin{equation*}
    \begin{aligned}
    \rho^S(G)z_r&= \sum_{\substack{t= 1 \\ t \neq r}}^{p}(2d_{rt}-1)z_t+ \sum_{\substack{t=p+1\\ t \neq r}}^{p+q}(2d_{rt}-1)z_t\\
    &\geq (p-1)3z_r+d_rz_s+5(q-d_r)z_s\\
    &\geq 3(p-1)z_r+(5q-4\Updelta_P)z_s.    \end{aligned}
\end{equation*}
That is, \begin{equation}\tag{9}
    (\rho^S(G)+3(1-p))z_r\geq (5q-4\Updelta_P)z_s.
\end{equation}\\
Similarly, for the component $z_s$ we have\\
\begin{equation*}
    \begin{aligned}
    \rho^S(G)z_s&= \sum_{\substack{t= 1 \\ t \neq s}}^{p}(2d_{st}-1)z_t+ \sum_{\substack{t=p+1\\ t \neq s}}^{p+q}(2d_{st}-1)z_t\\
    &\geq (q-1)3z_s+d_sz_r+5(p-d_s)z_r\\
    &\geq (q-1)3z_s+(5p-4\Updelta_Q)z_r.
    \end{aligned}
\end{equation*}

That is, \begin{equation}\tag{10}
    (\rho^S(G)+3(1-q))z_s\geq (5p-4\Updelta_Q)z_r.
\end{equation}\\
Combining (9) and (10), we get\\
$$(\rho^S(G)+3(1-p))(\rho^S(G)+3(1-q))\geq (5q-4\Updelta_P)(5p-4\Updelta_Q)$$
Since $z_t>0$ for $t= 1, \ldots, p+q$,\\
$$(\rho^S(G))^2+\rho^S(G)(6-3(p+q))+9(p-1)(q-1)-(5q-4\Updelta_P)(5p-4\Updelta_Q)\geq 0.$$
From this inequality, we get the result\\
$$\rho^S(G)\geq \frac{3n-6+\sqrt{9n^2-36pq+4(5p-4\Updelta_Q)(5q-4\Updelta_P)}}{2}.$$

For the equality case, we have $z_{r}=z_{t}$ for $t=1,2, \ldots, p$ and $z_{s}= z_{t}$ for $t=$ $p+1, p+2, \ldots, p+q$. That is, the eigenvector $z$ has maximum two different coordinates, $\Updelta_{P}$ is the degrees of vertices in $P$ and $\Updelta_{Q}$ is the degrees of vertices in $Q$, which implies that the graph $G$ is semi-regular.  If $G\neq K_{p,q}$, then $\Updelta_{P}<q$ and $\Updelta_{Q}<p$ (since $p \Updelta_{P}=q \Updelta_{Q}$) and the eccentricity $e(v_i)= 3$ for every $v_i\in V(G).$  
\end{proof}

\section{Distance Seidel energy and its bounds}
\begin{proposition}
The distance Seidel energy of some standard graphs are given by,
\begin{itemize}
    \item [(i)] $E_{\mathcal{D}^S}(K_n)= 2n-2= E_{\mathcal{D}}(K_n)= E(K_n),$
    \item[(ii)] $E_{\mathcal{D}^S}(K_n-e)= 2n$
    \item[(iii)]For $a, b>1$, $E_{\mathcal{D}^S}(K_{a, b})= 6(a+b-2),$
    \item[(iv)] $E_{\mathcal{D}^S}(K_{n, n})= 6(2n-2)= 6E_{\mathcal{D}^S}(K_n),$
    \item[(v)] $E_{\mathcal{D}^S}(S_n)= 3n-6+\sqrt{9n^2-32n+32}.$
\end{itemize}
\end{proposition}
\begin{corollary}
For $n>1$, $E_{\mathcal{D}^S}(K_n-e)> E_{\mathcal{D}^S}(K_n)$
\end{corollary}
\begin{theorem}
For $a, b>1$, the distance Seidel characteristic polynomial of $K_{a,b}-e$, where $e$ is any edge of $K_{a,b}$, is given by\\
$p(x)= (-1)^{a+b}(x-3)^{a+b-4}(-x^4+x^3(12-3a-3b)+x^2(-30+27a+27b-8ab)+x(-132-33a-33b+48ab)+(551-191a-191b+56ab).$
\end{theorem}
\begin{proof}
Consider partition the vertex set of $K_{a, b}$ as the union of the sets $\{u_1,\ldots, u_a\}$ and $\{u_{a+1}, \ldots, u_{a+b}\}$. Also let, $e= u_au_{a+1}$
be an edge of $K_{a, b}$. Then the distance Seidel matrix of $K_{a, b}-e$ has the form\\

$\mathcal{D}^S(K_{a, b}-e)= \left(\begin{array}{cccccccccc}
0 & -3 & \cdots & -3 & -3 & -1 & -1 & \cdots & -1 & -1 \\
-3 & 0 & \cdots & -3 & -3 & -1 & -1 & \cdots & -1 & -1 \\
\vdots & \vdots & \ddots & \vdots & \vdots & \vdots & \vdots &\ddots & \vdots & \vdots \\
-3 & -3 & \cdots & 0 & -3 & -1 & -1 & \cdots & -1 & -1 \\
-3 & -3 & \cdots & -3 & 0 & -5 & -1 & \cdots & -1 & -1 \\
-1 & -1 & \cdots & -1 & -5 & 0 & -3 & \cdots & -3 & -3 \\
-1 & -1 & \cdots & -1 & -1 & -3 & 0 & \cdots & -3 & -3 \\
\vdots & \vdots & \ddots& \vdots & \vdots & \vdots & \vdots & \ddots& \vdots & \vdots \\
-1 & -1 & \cdots & -1 & -1 & -3 & -3 & \cdots & 0 & -3 \\
-1 & -1 & \cdots & -1 & -1 & -3 & -3 & \cdots & -3 & 0
\end{array}\right)_{(a+b)\times (a+b)}.$\\

Then the characteristic polynomial of $\mathcal{D}^S(K_{a, b}-e)-3I$ is,
\\
$\begin{aligned}
&\operatorname{det}\left(xI-\mathcal{D}^S\left(K_{a, b}-e\right)+3I\right)\\ 
&= \operatorname{det}\left(\begin{array}{cccc|cccc}
x+3 & 3 & \cdots & 3 & 1 & 1 & \cdots & 1 \\
3 & x+3 & \cdots & 3 & 1 & 1 & \cdots & 1 \\
\vdots & \vdots & \ddots& \vdots & \vdots & \vdots & \ddots& \vdots \\
3& 3 & \cdots & x+3 & 5 & 1 & \cdots & 1 \\
\hline 1 & \cdots & 1 & 5 & x+3 & 3 & \cdots & 3 \\
1 & \cdots & 1 & 1 & 3 & x+3 & \cdots & 3 \\
\vdots &\ddots & \vdots & \vdots & \vdots & \vdots & \ddots & \vdots \\
1 & \cdots & 1 & 1 & 3 & 3 & \cdots & x+3
\end{array}\right)_{(a+b)\times (a+b)}
\end{aligned}$\\
\\
$\begin{aligned}
&= \operatorname{det}\left(\begin{array}{cccc|cccc}
x & 0 & \cdots & 0 & -8 & -8 & \cdots & -3x-8 \\
0 & x & \cdots & 0 & -8 & -8 & \cdots & -3x-8 \\
\vdots & \vdots & \ddots & \vdots & \vdots & \vdots & \ddots & \vdots \\
0& 0 & \cdots & x & -4 & -8 & \cdots & -3x-8 \\
\hline 0 & \cdots & 0 & 4 & x & 0 & \cdots & -x \\
0 & \cdots & 0 & 0 & 0 & x & \cdots & -x \\
\vdots & \ddots & \vdots & \vdots & \vdots & \vdots & \ddots & \vdots \\
1 & \cdots & 1 & 1 & 3 & 3 & \cdots & x+3
\end{array}\right)_{(a+b)\times (a+b)}
\end{aligned}$\\
\\
\\
{\scriptsize$\begin{aligned}
&= (-1)^{a+b}\operatorname{det}\left(\begin{array}{cccc|ccc|ccc}
-x & -x & \cdots & -x & -x & -8-3x & -8-3x & \cdots & -8-3x & -x^{2}-6x-8 \\
x & 0 & \cdots & 0 & 0 & -8 & -8 & \cdots & -8 & -8-3x \\
\vdots & \vdots &\ddots & \vdots & \vdots & \vdots & \vdots & \ddots & \vdots & \vdots \\
0 & 0 & \cdots & 0 & 0 & -8 & -8 & \cdots & -8 & -8-3x \\
\hline 0 & 0 & \cdots & x & 0 & -8 & -8 & \cdots & -8 & -8-3x \\
0 & 0 & \cdots & 0 & x & -4 & -8 & \cdots & -8 & -8-3x \\
0 & 0 & \cdots & 0 & 4 & x & 0 & \cdots & 0 & -x \\
\hline 0 & 0 & \cdots & 0 & 0 & 0 & -x & \cdots & 0 & x \\
\vdots & \vdots & \ddots& \vdots & \vdots & \vdots & \vdots & \ddots& \vdots & \vdots \\
0 & 0 & \cdots & 0 & 0 & 0 & 0 & \cdots & x & -x 
\end{array}\right)_{(a+b-1)\times (a+b-1)}
\end{aligned}$}\\
\\
\\
$\begin{aligned}
&=(-1)^{a+b} x^{a-2} x^{b-2} \operatorname{det}\left(\begin{array}{ccc}
-x & r & s \\
x & -4 & t \\
4 & x & -x
\end{array}\right)
\end{aligned}$, where\\

$$r= -8-3x+(-8)(a-2),$$
$$s= (-x^2-6x-8)+(-8)(a-2)+(b-2)(-8-3x+(-8)(a-2)),$$
$$t= -8-3x+(-8)(b-2).$$
Then 
$det\left(xI-\mathcal{D}^S\left(K_{a, b}-e\right)+3I\right)
= (-1)^{a+b}x^{a+b-4}(-x^4-3x^3(a+b)+8x^2(3-ab)+48x(a+b-2)+128(a-1)(b-1)).$\\
Now shift the eigenvalues by $3$ to obtain the characteristic polynomial for $\mathcal{D}^S(K_{a, b}-e)$,
$det(\mathcal{D}^S(K_{a, b}-e)-xI)= (-1)^{a+b}(x-3)^{a+b-4}(-x^4+x^3(12-3a-3b)+x^2(-30+27a+27b-8ab)+x(-132-33a-33b+48ab)+(551-191a-191b+56ab)).$
\end{proof}

\begin{proposition}
For $n\geq 2$, $E_{\mathcal{D}^S}(K_{2,n}-e)> E_{\mathcal{D}^S}(K_{2,n}).$

\end{proposition}
\begin{proof}
From Theorem $5.1$, it follows that\\
$\operatorname{det}\left(xI-\mathcal{D}^S(K_{2,n}-e)\right)= (-1)^n(x-3)^{n-2}(-x^4+x^3(6-3n)+x^2(24+11n)+x(-198+63n)+169-79n).$\\
Let $\beta_4\leq \beta_3\leq \beta_2\leq \beta_1$ be the four roots of the polynomial,
\begin{equation}\tag{11}
s(x)= -x^4+x^3(6-3n)+x^2(24+11n)-x(198-63n)+169-79n.
\end{equation}
\\
For $n=2,3$ we have,
\begin{equation*}
    \begin{aligned}
    E_{\mathcal{D}^S}(K_{2,2}-e)&= 14.94> 8= E_{\mathcal{D}^S}(K_{2,2}).\\
    E_{\mathcal{D}^S}(K_{2,3}-e)&= 20.41> 19.627= E_{\mathcal{D}^S}(K_{2,3}).
    \end{aligned}
\end{equation*}
Therefore, let us consider the case when $n\geq 4.$\\
From (11) we have, $\beta_1\beta_2\beta_3\beta_4= -169+79n> 0.$ Then the number of positive roots of $s(x)$ can be $0, 2$ or $4.$ Also, $\sum \beta_i\beta_j= -24-11n< 0,$ that is, all $\beta_i$'s can't be positive. Moreover, $\beta_1+\beta_2+\beta_3+\beta_4= 6-3n< 0.$ That is all $\beta_i$'s can't be negative.\\
Therefore, $s(x)$ must possess exactly two positive and two negative roots. Now $s(1)= -8n<0$ and $s(2)= 67n-99> 0$ since $n\geq 4.$ Then by intermediate value theorem we get, $1<\beta_2< 2.$\\
Similarly, $s(6)= 47n-155>0$ and $s(7)= -384-128n<0$, implies $6<\beta_1<7.$\\
Now consider,\\
\begin{equation*}
    \begin{aligned}
    E_{\mathcal{D}^S}(K_{2,n}-e)&= 3(n-2)+|\beta_1|+|\beta_2|+|\beta_3|+|\beta_4|\\
&= 6n-12+2(\beta_1+\beta_2)\\
&>6n-12+14= 6n+2.
    \end{aligned}
\end{equation*}
Then $E_{\mathcal{D}^S}(K_{2,n}-e)> E_{\mathcal{D}^S}(K_{2,n})$.

\end{proof}
\begin{proposition}
For $n\geq 2$, $E_{\mathcal{D}^S}(K_{3,n}-e)> E_{\mathcal{D}^S}(K_{3,n}).$

\end{proposition}
\begin{proof}
From Theorem $5.1$, it follows that\\
$\operatorname{det}\left(xI-\mathcal{D}^S(K_{3,n}-e)\right)= (-1)^{n+1}(x-3)^{n-1}(-x^4+x^3(3-3n)+x^2(51+3n)+x(-231+111n)-22-23n).$\\
Let $\beta_4\leq \beta_3\leq \beta_2\leq \beta_1$ be the four roots of the polynomial
\begin{equation}\tag{12}
    s(x)= -x^4+x^3(3-3n)+x^2(51+3n)+x(-231+111n)-22-23n.
\end{equation}
\\
For $n=2,3$ we have,
\begin{equation*}
    \begin{aligned}
    E_{\mathcal{D}^S}(K_{3,2}-e)&= 20.41> 19.627= E_{\mathcal{D}^S}(K_{3,2}).\\
    E_{\mathcal{D}^S}(K_{3,3}-e)&= 25.6> 24= E_{\mathcal{D}^S}(K_{3,3}).
    \end{aligned}
\end{equation*}

Consider the case when $n\geq 4.$\\
From (12) we have,\\ 
$$\beta_1\beta_2\beta_3\beta_4= 22+23n> 0.$$
Then the number of positive roots of $s(x)$ can be $0, 2\;  \text{or}\;  4.$\\
Now by similar arguments as in Proposition $5.2$, we obtain
\begin{equation*}
\begin{aligned}
E_{\mathcal{D}^S}(K_{3,n}-e)&= 3(n-1)+|\beta_1|+|\beta_2|+|\beta_3|+|\beta_4|\\
&= 6n-6+2(\beta_1+\beta_2)\\
&>6n+6= E_{\mathcal{D}^S}(K_{3,n}).
\end{aligned}
\end{equation*}

\end{proof}
\begin{problem}
Let $G= K_{a,b}, \;\text{with}\; a, b> 1$. Is $E_{\mathcal{D}^S}(K_{a,b}-e)>E_{\mathcal{D}^S}(K_{a,b})$, for any edge $e$ of $K_{a,b}?$ 
\end{problem}

\begin{theorem}
Let $G$ be a connected graph of order $n$. Then for every $1\leq r\leq n$, $|\partial_r^S|\leq \frac{E_{\mathcal{D}^S}(G)}{2}.$ Moreover if $G\cong K_n$, then $|\partial_n^S|= \frac{E_{\mathcal{D}^S}(G)}{2}$, where $\partial_n^S$ is the smallest eigenvalue of $\mathcal{D}^S(G).$
\end{theorem}
\begin{proof}
Since the diagonal entries of $\mathcal{D}^S(G)$ are all zeros, $\sum_{r=1}^{n}\partial_r^S= 0.$\\
\begin{equation*}
    \begin{aligned}
    &\partial_t^S= -\sum_{r=1, r\neq t}^{n}\partial_r^S\\
    &|\partial_t^S|= \Big| \sum_{r=1, r\neq t}^{n}\partial_r^S\Big|\leq \sum_{r=1, r\neq t}^{n}|\partial_r^S|
    \end{aligned}
\end{equation*}

adding $|\partial_t^S|$ on both sides, we get\\
$$2|\partial_t^S|\leq \sum_{r=1, r\neq t}^{n}|\partial_r^S|+|\partial_t^S|= E_{\mathcal{D}^S}(G)$$\\

If $G\cong K_n$, then clearly $|\partial_n^S|= \frac{E_{\mathcal{D}^S}(G)}{2}$, since $|\partial_n^S|=n-1.$ 
\end{proof}
\begin{theorem}
Let $G$ be a $k$-transmission regular graph of order n with number of non negative and negative $\mathcal{D}$-eigenvalues as $a^{+}$ and $a^{-}$ respectively. Suppose $\partial_1=k, \partial_2, \ldots, \partial_n$ and $n-1-2k, -1-2\partial_2, \ldots, -1-2\partial_n$ be the $\mathcal{D}$-eigenvalues and $\mathcal{D}^S$-eigenvalues of $G$ respectively. If $\partial_r\notin (-1, 0) \; \text{for} \; 2\leq r \leq n$, then $$E_{\mathcal{D}^S}(G)= 2(a^{+}-n+E_{\mathcal{D}}(G)).$$
\end{theorem}
\begin{proof}
Suppose that $\partial_{n}\leq \partial_{n-1}\leq \cdots \leq \partial_{1}= k.$
Clearly, $k\geq \frac{(n-1)}{2}$, that is, $n-1-2k<0.$ \\
For $r= a^{+}+1,\ldots, n$, we have $\partial_r\leq -1$, which implies $-1-2\partial_r\geq 1.$\\

Therefore, the non-negative distance Seidel eigenvalues are $\partial_r^S= -1-2\partial_r$ for $r= a^{+}+1, \ldots, n$. Then
\begin{equation*}
    \begin{aligned}
    E_{\mathcal{D}^S}(G)&= 2\sum_{r= a^{+}+1}^{n}\partial_r^S\\
    &= 2\sum_{r= a^{+}+1}^{n}-1-2\partial_r\\
    &= -2a^{-}+2E_{\mathcal{D}}(G).
    \end{aligned}
\end{equation*}

Now by replacing $a^-$ with $n-a^{+}$, we get
$$E_{\mathcal{D}^S}(G)= 2(a^{+}-n+E_{\mathcal{D}}(G)).$$

\end{proof}
\begin{theorem}
Let $G$ be a connected graph of order n, with number of non negative and negative $\mathcal{D}$-eigenvalues as $a^{+}$ and $a^{-}$ respectively. If $\partial_r\notin (-1 / 2,0)\; \text{for}\; r= a^{+}+1, \ldots, n$ and $\partial_{r} \geq\frac{(n-1)}{2}, \; \text{for}\; r= 1, \ldots, a^{+}$, then 
$$
E_{\mathcal{D}^S}(G) \geq 2 E_{\mathcal{D}}(G)-2 a^{-}.
$$
\end{theorem}
\begin{proof}
If $\partial_1, \ldots, \partial_n$ and $\partial_{1}^{S}, \ldots, \partial_{n}^{S}$ be the $\mathcal{D}$-eigenvalues and $\mathcal{D}^S$-eigenvalues respectively. Let $\partial_r\notin (-1 / 2,0)$ for $r= a^{+}+1, \ldots, n$ and $\partial_{r} \geq\frac{(n-1)}{2}, r= 1, \ldots, a^{+}$, then we have $\partial_{1} \geq \cdots \geq$ $\partial_{a^{+}} \geq\frac{(n-1)}{2}$ and $-1-2 \partial_{r} \geq 0$ for $r= a^{+}+1, \ldots, n$. From Theorem $3.2$, we have $\partial^{S}_{n-r+1} \geq-1-2 \partial_{r} \geq 0$ for $r= a^{+}+1, \ldots, n$. This implies that
$$
\sum_{a^{+}+1}^{n} \partial^{S}_{n-r+1} \geq \sum_{a^{+}+1}^{n}\left(-1-2 \partial_{r}\right)
$$
That is,

\begin{equation*}
    \begin{aligned}
    \frac{E_{\mathcal{D}^S}(G)}{2} &\geq-\sum_{a^{+}+1}^{n} 1-2 \sum_{a^{+}+1}^{n} \partial_{n-r+1}\\
    &\geq-a^{-}+E_{\mathcal{D}}(G).
    \end{aligned}
\end{equation*}
Hence $E_{\mathcal{D}^S}(G) \geq 2 E_{\mathcal{D}}(G)-2 a^{-}.$
\end{proof}
In a  connected graph $G$ of order $n\geq 2$, we have $E_{\mathcal{D}^S}(G)\geq 2\rho^S (G),$ and equality holds if and only if $G$ has exactly one negative $\mathcal{D}^S$-eigenvalue. Then the lower bounds for $\rho^S (G)$ may be converted to lower bounds for $E_{\mathcal{D}^S}(G).$
\begin{theorem}
For a connected graph $G$ of order $n$,\\
$$E_{\mathcal{D}^S}(G)\geq 2\sqrt{\frac{\sum_{r=1}^{n}(\mathcal{D}^S_{r})^2}{n}}$$
where $\mathcal{D}^S_r= \sum_{\substack{t=1\\ t\neq r}}^{n}2d_{rt}-1$ for $r= 1, 2, \ldots, n$. Equality holds if and only if $\mathcal{D}^S$ has a constant row sum and $G$ has exactly one negative $\mathcal{D}^S$-eigenvalue.
\end{theorem}
\begin{proof}
The proof follows from Theorem $4.2$.
\end{proof}
\begin{proposition}
For a connected graph $G$ of order $n$ with maximum degree $\Delta$ and second maximum degree $\Delta'$,\\
$$E_{\mathcal{D}^S}(G)\geq 2\sqrt{(3n-2\Delta-3)(3n-2\Delta'-3)}.$$
\end{proposition}
\begin{proof}
The proof follows from Theorem $4.3.$
\end{proof}

\begin{theorem}
Let $G$ be a connected graph of order $n$ and  size $m$. Then\\
$$\sqrt{2T+n(n-1)|det(\mathcal{D}^S(G))|^{2/n}}\leq E_{\mathcal{D}^S}(G)\leq \sqrt{2nT},\; \text{where}\; T= \sum_{r<t}(1-2d_{rt})^2 $$
\end{theorem}
\begin{proof}
Consider \begin{equation*}
    \begin{aligned}
    (E_{\mathcal{D}^S}(G))^2&= (\sum_{r=1}^{n}|\partial_r^S|)^2\\
    &= \sum_{r=1}^{n}(\partial_r^S)^2+ 2\sum_{r<t}|\partial_r^S||\partial_t^S|\\
    &= 2T+2\sum_{r<t}|\partial_r^S||\partial_t^S|\\
    &= 2T+\sum_{r\neq t}|\partial_r^S||\partial_t^S|.
    \end{aligned}
\end{equation*}
Using the AM-GM inequality, we get
\begin{equation*}
    \begin{aligned}
    \frac{1}{n(n-1)}\sum_{r\neq t}|\partial_r^S||\partial_t^S|&\geq (\prod_{r\neq t}|\partial_r^S||\partial_t^S|)^{\frac{1}{n(n-1)}}\\
    &= (\prod_{r=1}^{n}|\partial_r^S|^{2(n-1)})^{\frac{1}{n(n-1)}}\\
    &= |det(\mathcal{D}^S(G))|^{2/n}.
    \end{aligned}
\end{equation*}
That is, $$\sum_{r\neq t}|\partial_r^S||\partial_t^S|\geq n(n-1)|det(\mathcal{D}^S(G))|^{2/n}.$$
Then \begin{equation*}
    \begin{aligned}
    (E_{\mathcal{D}^S}(G))^2&= 2T+\sum_{r\neq t}|\partial_r^S||\partial_t^S|\\
    &\geq 2T+n(n-1)|det(\mathcal{D}^S(G))|^{2/n}.
    \end{aligned}
\end{equation*}
Now consider,
\begin{equation*}
    \begin{aligned}
    X&= \sum_{r=1}^{n}\sum_{t=1}^{n}(|\partial_r^S|-|\partial_t^S|)^2\\
    &= n\sum_{r=1}^{n}(\partial_r^S)^2+n\sum_{t=1}^{n}(\partial_t^S)^2-2(\sum_{r=1}^{n}|\partial_t^S|)(\sum_{t=1}^{n}|\partial_t^S|)\\
    &= 4nT-2E_{\mathcal{D}^S}(G)^2.
    \end{aligned}
\end{equation*}

Then $X\geq 0$ implies that $4nT-2E_{\mathcal{D}^S}(G)^2\geq 0,$\\
$$E_{\mathcal{D}^S}(G)\leq \sqrt{2nT}.$$
Combining two inequalities we get,
$$\sqrt{2T+n(n-1)|det(\mathcal{D}^S(G))|^{2/n}}\leq E_{\mathcal{D}^S}(G)\leq \sqrt{2nT}.$$

\end{proof}
Let $G$ be a graph of order $n$ and size $m$ with $\tilde{\omega}\leq 2$. Then in $G,\; m$ pairs of vertices are at a distance $1$, and $nC_2-m$ pairs of vertices are at a distance $2.$ Then we have $T= \frac{9n^2-9n-16m}{2}.$
Now by using Theorem $5.12$ we get the following corollary.
\begin{corollary}
For a connected graph $G$ of order $n$ and size $m$ with $\tilde{\omega}\leq 2$,
{\small$$\sqrt{9n^2-9n-16m+n(n-1)|det(\mathcal{D}^S(G))|^{2/n}}\leq E_{\mathcal{D}^S}(G)\leq \sqrt{9n^3-9n^2-16mn}.$$}
\end{corollary}

\section{Distance Seidel spectrum of some graph operations}
In this section we obtain the $\mathcal{D}^S$-spectrum of graphs formed by using different graph operations.
Consider two graphs $G_1$ and $G_2$ having disjoint vertex sets $V(G_1)$ and $V(G_2)$, then the join of $G_1$ and $G_2$, denoted by $G_1\vee G_2$, is the graph obtained from $G_1\cup G_2$ by adding an edge joining each vertex of $G_1$ to every vertex of $G_2$.
\begin{theorem}
For $r=1, 2$, let $G_r$ be an $k_r$-regular graph of order $n_r$ and size $m_r$ with adjacency eigenvalues $\lambda_{r,1}= k_r, \lambda_{r,2}, \ldots, \lambda_{r,{n_r}}$. Then the $\mathcal{D}^S$-eigenvalues of $G_1\vee G_2$ are given by
\begin{itemize}
    \item [(i)] $3+2\lambda_{1,t}$, for $t= 2, 3, \ldots, n_1$
    \item [(ii)] $3+2\lambda_{2,t}$, for $t= 2, 3, \ldots, n_2$
    \item [(iii)] $\frac{-(3n_1+3n_2-2k_1-2k_2-6)\pm \sqrt{(3n_1-3n_2-2k_1+2k_2)^2+4n_1n_2}}{2}.$
\end{itemize}

\end{theorem}
\begin{proof}
The distance Seidel matrix of $G_1\vee G_2$ is
$$\mathcal{D}^S(G_1\vee G_2)= \begin{bmatrix}3(I-J)+2A(G_1)& -J\\ -J& 3(I-J)+2A(G_2)\end{bmatrix}.$$
Clearly $k_1$ is the Perron eigenvalue of $A(G_1)$ (since $G_1$ is a $k_1$-regular graph). Now let $\lambda\neq k_1$ be an eigenvalue of $A(G_1)$ with an eigenvector U. Then by Perron-Frobenius theorem $JU= 0.$\\
Let $\phi= \begin{bmatrix}U\\0\end{bmatrix}$, then $$\mathcal{D}^S\phi= (3+2\lambda) \phi.$$
Therefore, for every adjacency eigenvalue $\lambda$ of $G_1$, we get a distance Seidel eigenvalue $3+2\lambda$ of $G_1\vee G_2$ with an eigenvector $\phi.$ \\
Similarly for every adjacency eigenvalue $\mu$ of $G_2$, we get a distance Seidel eigenvalue $3+2\mu$ of $G_1\vee G_2$ with an eigenvector $\begin{bmatrix}0\\V\end{bmatrix}$ where $A(G_2)V= \lambda_{2,r} V$ (for $r= 2,3,...,n_2$).\\
In total we get $n_1+n_2-2$ eigenvalues of $G_1\vee G_2.$ Now consider the equitable quotient matrix $M$ of $\mathcal{D}^S(G_1\vee G_2)$. Then the rest of the eigenvalues are the eigenvalues of the matrix
$$M= \begin{bmatrix}3(1-n_1)+2k_1& -n_2\\-n_1& 3(1-n_2)+2k_2\end{bmatrix}.$$\\
$|xI-M|= x^2+x((3n_1-2k_1-3)+(3n_2-2k_2-3))+(3n_1-2k_1-3)(3n_2-2k_2-3)-n_1n_2.$\\
Then $x= \frac{-(3n_1+3n_2-2k_1-2k_2-6)\pm \sqrt{(3n_1-3n_2-2k_1+2k_2)^2+4n_1n_2}}{2}.$
\end{proof}
Next corollary is an application of the above theorem.
\begin{corollary}
The distance Seidel spectrum of the wheel graph $W_n$, the complete split graph $CS(n,p)$, the friendship graph $f_n$, and $  K_{\underbrace{n, n, \ldots, n}_{\substack{q-times}}}\; \left(q\leq n\right)$  are given by,

{\small\begin{itemize}
\item [(i)]$Spec_{\mathcal{D}^S}(W_n)= \begin{pmatrix}3+4\cos{\frac{2\pi t}{n-1}}& \frac{10-3n\pm \sqrt{9n^2-56n+96}}{2}\\1& 1\end{pmatrix},\; \text{for}\; 1\leq t \leq n-2,$
\item[(ii)]$Spec_{\mathcal{D}^S}(CS(n,p))= \begin{pmatrix}1& 3& \frac{2p-3n+4\pm \sqrt{12p^2+9n^2+16p-12n-20np+4}}{2}\\p-1& n-p-1& 1\end{pmatrix},$
\item[(iii)]$Spec_{\mathcal{D}^S}(f_n)= \begin{pmatrix}1& 5& \frac{5-6n\pm \sqrt{36n^2-52n+25}}{2}\\n& n-1& 1\end{pmatrix},$
\item [(iv)]$Spec_{\mathcal{D}^S}(  K_{\underbrace{n, n, \ldots, n}_{\substack{q-times}}} )= \begin{pmatrix}3& 3-2n& 3-n(q+2)\\ q(n-1)& q-1& 1\end{pmatrix},$ where $q\leq n.$\\
    In particular, the $\mathcal{D}^S$-spectrum of the cocktail party graph $CP(n)$ is given by,\\
    $Spec_{\mathcal{D}^S}(  K_{\underbrace{2, 2, \ldots, 2}_{\substack{n-times}}} )= \begin{pmatrix}3& -1& -1-2n\\ n& n-1& 1\end{pmatrix}.$
\end{itemize}}
\end{corollary}
\begin{proof}
We have $W_n= C_{n-1}\vee K_1$, $CS(n,p)= K_p\vee \Bar{K}_{n-p}$, $f_n= K_1\vee nK_2$, and $K_{\underbrace{n, n, \ldots, n}_{\substack{q-times}}}= K_{\underbrace{n, n, \ldots, n}_{\substack{(q-1)-times}}}\vee \Bar{K_n}$ then the proof follows from Theorem $6.1.$
\end{proof}
\begin{corollary}
The distance Seidel energy of the complete split graph $CS(n,p)$, the friendship graph $f_n$, and  $  K_{\underbrace{n, n, \ldots, n}_{\substack{q-times}}} $  ($q\leq n$) are given by
\begin{itemize}
\item[(i)] $E_{\mathcal{D}^S}(CS(n,p))= 2p-3n+4-\sqrt{12p^2+9n^2+16p-12n-20np+4},$
 \item [(ii)] $E_{\mathcal{D}^S}(f_n)= 6n-5+\sqrt{36n^2-52n+25},$
    \item [(iii)] $E_{\mathcal{D}^S}(K_{\underbrace{n, n, \ldots, n}_{\substack{q-times}}})= 3q(n-1),$ where $q\leq n.$
    
\end{itemize}
\end{corollary}
Theorem $6.1$ gives the $\mathcal{D}^S$- spectrum of join of two graphs when both $G_1$ and $G_2$ are regular. Through the next theorem we give the $\mathcal{D}^S$- spectrum of $G_0$ over union of two graphs $G_1$ and $G_2$, in which $G_1\cup G_2$ can be non regular. 
\begin{theorem}
For $r= 0, 1, 2$, let $G_r$ be a $k_r$-regular graph of order $n_r$ with adjacency eigenvalues $k_r=\lambda_{r1}\geq \lambda_{r2}\geq \ldots \geq \lambda_{rn_{r}}$. Then $\mathcal{D}^S$-eigenvalues of $G_0\vee (G_1\cup G_2)$ are given by,
\begin{itemize}
    \item [(i)] $3+2\lambda_{rt}, t= 2,3,\ldots, n_r, r= 0, 1, 2.$
    \item [(ii)]the eigenvalues of the matrix
    \[\begin{bmatrix}3-3n_0+2k_0& -n_1& -n_2\\ -n_0& 3-3n_1+2k_1& -3n_2\\ -n_0& -3n_1& 3-3n_2+2k_2\end{bmatrix}.\]
\end{itemize}
\end{theorem}
\begin{proof}
The distance Seidel matrix of $G_0\vee (G_1\cup G_2)$ is\\
{\scriptsize$$\mathcal{D}^S(G_0\vee (G_1\cup G_2))= \begin{bmatrix}3(I-J)+2A(G_0)& -J& -J\\ -J& 3(I-J)+2A(G_1)& -3J\\ -J& -3J& 3(I-J)+2A(G_2)\end{bmatrix}.$$}\\
By similar arguments in Theorem $6.1$ , corresponding to the adjacency eigenvalues $\lambda_{rt}$ of $A(G_r)$, we get $\mathcal{D}^S$-eigenvalues $3+2\lambda_{rt}\;(t= 2,3,\ldots, n_r, r= 0, 1, 2).$\\
That is in total we get $n_1+n_2+n_3-3$ eigenvalues of $\mathcal{D}^S(G_0\vee (G_1\cup G_2)).$ Now consider the equitable quotient matrix $M$ of $\mathcal{D}^S(G_0\vee (G_1\cup G_2)).$ Then the rest of the eigenvalues are the eigenvalues of the matrix\\
\[M=\begin{bmatrix}3-3n_0+2r_0& -n_1& -n_2\\ -n_0& 3-3n_1+2r_1& -3n_2\\ -n_0& -3n_1& 3-3n_2+2r_2\end{bmatrix}.\]
\end{proof}

 Let $G$ be a graph with $V(G) = \{v_1, v_2,\ldots ,v_n\}$. Make another copy of $G$ with vertices $\{v_{1}^{'}, v_{2}^{'},\ldots ,v_{n}^{'}\}$ in which $v_{r}^{'}$ corresponds to $v_r$ for each $r$ in such a way that for each $r, \; v_{r}^{'}$ is adjacent to every vertices in the neighborhood of $v_r$ in $G$. The resultant graph is called the \textit{double graph} \cite{indulal2006pair} of $G$, denoted by $D_2G.$

\begin{theorem}
Let $G$ be a graph of order $n$ with $\mathcal{D}^S$-eigenvalues $\partial_1^{S}, \partial_2^{S}, \ldots, \partial_n^{S}.$ Then
$$Spec_{D^{S}}(D_2 G)= \begin{pmatrix}2\partial_r^{S}-3& 3\\1& n\end{pmatrix}, \; 1\leq r\leq n.$$
\end{theorem}
\begin{proof}
The distance Seidel matrix of $D_2 G$ is\\
$$\mathcal{D}^S(D_2 G)= \begin{bmatrix}\mathcal{D}^S& \mathcal{D}^S-3I\\ \mathcal{D}^S-3I& \mathcal{D}^S\end{bmatrix}.$$
Then the proof follows from Lemma $2.1$.
\end{proof}
\begin{corollary}
For a graph $G$ of order $n$, $E_{\mathcal{D}^S}(D_2G)\leq 2E_{\mathcal{D}^S}(G)+6n.$
\end{corollary}
For two graphs $G_1 = (V(G_1), E(G_1))$ and $G_2 = (V(G_2), E(G_2))$, their \textit{cartesian product} $G_1\times G_2$ is a graph with $V(G_1\times G_2)= V(G_1)\times V(G_2)$ in which two vertices ($v_1, u_1$) and ($v_2, u_2$) are adjacent if $v_1= v_2$ and ($u_1, u_2$) $\in  E(G_2)$ or $u_1 = u_2$ and
($v_1, v_2$) $\in E(G_1)$.
\begin{theorem}
Let $G$ be a distance regular graph with $\mathcal{D}$-eigenvalues $\partial_1, \ldots, \partial_n.$ Then
$$Spec_{D^{S}}(G\times K_2)= \begin{pmatrix}-1-4\partial_r& 2n-1& -1\\1& 1& n-1\end{pmatrix}, \; 1\leq r \leq n.$$

\end{theorem}
\begin{proof}
The distance Seidel matrix of $G\times K_2$ is\\
$$\mathcal{D}^S(G\times K_2)= \begin{bmatrix}J-I-2\mathcal{D}& -J-2\mathcal{D}\\ -J-2\mathcal{D}& J-I-2\mathcal{D}\end{bmatrix}.$$
Then the proof follows from Lemma $2.1$.
\end{proof}
\begin{corollary}
Let $G$ be a distance regular graph, then $E_{\mathcal{D}^S}(G\times K_2)\leq 4E_{\mathcal{D}}(G)+4n-2.$
\end{corollary}
 For two graphs $G_1 = (V(G_1), E(G_1))$ and $G_2 = (V(G_2), E(G_2))$, their \textit{lexicographic product} $G_1[G_2]$ is the graph with $V(G_1[G_2])= V(G_1)\times V(G_2)$ in which two vertices ($v_1, u_1$) and ($v_2, u_2$) are adjacent if ($v_1, v_2$) $\in  E(G_1)$ or $v_1 = v_2$ and
($u_1, u_2$) $\in E(G_2)$.
\begin{theorem}
Let $G$ be a graph of order $n$ with $\mathcal{D}^S$-eigenvalues $\partial_1^{S}, \partial_2^{S}, \ldots, \partial_n^{S}.$ Then
$$Spec_{D^{S}}(G[K_2])= \begin{pmatrix}2\partial_r^{S}-1& 1\\1& n\end{pmatrix}, \; 1\leq r \leq n.$$
\end{theorem}
\begin{proof}
The distance Seidel matrix of $G[K_2]$ is\\
$$\mathcal{D}^S(G[K_2])= \begin{bmatrix}\mathcal{D}^S& \mathcal{D}^S-I\\ \mathcal{D}^S-I& \mathcal{D}^S\end{bmatrix}.$$
Then the proof follows from Lemma $2.1$.

\end{proof}
\begin{corollary}
For a graph $G$ of order $n$, $E_{\mathcal{D}^S}(G[K_2])\leq 2E_{\mathcal{D}^S}(G)+2n.$
\end{corollary}
Let $G$ be a graph with $V(G)= \{v_1, v_2,\ldots ,v_n\}$. Then the \textit{extended double cover graph} (EDC-graph)  \cite{alon1986eigenvalues} of $G$ is a bipartite graph $X$ defined as follows, $V(X)= \{v_1, v_2,\ldots , v_n, u_1, u_2,\ldots , u_n\}$ in which $v_r$ is adjacent to $u_r$ for each
$r = 1, 2,\ldots, n$ and $v_r$ is adjacent to $u_t$ if $v_r$ is adjacent to $v_t$ in $G$ .
\begin{theorem}
Let $G$ be a $k$-regular graph of order $n$ and diameter $2$ with $A$-eigenvalues $k=\lambda_1\geq \lambda_2 \geq \cdots \geq \lambda_n.$ Then \\
$$Spec_{D^{S}}(EDC(G))= \begin{pmatrix}-8n+4k+7& 7+4\lambda_r& 2n-4k-1& -1-4\lambda_r\\1& 1& 1& 1\end{pmatrix},$$ where $r= 2, \ldots, n.$
\end{theorem}
\begin{proof}
The distance Seidel matrix of $EDC(G)$ is\\
$$\mathcal{D}^S(EDC(G)= \begin{bmatrix}3I-3J& 4I-5J+4A\\ 4I-5J+4A& 3I-3J\end{bmatrix}.$$
Then the proof follows from Lemma $2.1$.
\end{proof}

\section{Distance Seidel cospectral and integral graphs}
Two graphs are said to be \textit{distance Seidel cospectral} ($\mathcal{D}^S$-cospectral) if they have the same set of $\mathcal{D}^S$-eigenvalues.
\begin{proposition}
Two regular graphs with the same degree and diameter $2$ have cospectral distance Seidel matrices if and only if they have cospectral adjacency matrices.
\end{proposition}
The following Corollary's give different families of $\mathcal{D}^S$-cospectral graphs.

\begin{corollary}
Let $G$ and $G^{'}$ be two $D$-cospectral transmission regular graphs with same transmission regularity. Then $G$ and $G^{'}$ are $\mathcal{D}^S$-cospectral.
\end{corollary}
\begin{corollary}
Let $G$ and $G^{'}$ be two $D$-cospectral graphs, then $G\times K_2$ and $G^{'}\times K_2$ are $\mathcal{D}^S$-cospectral.
\end{corollary}
\begin{corollary}
Let $G$ and $G^{'}$ be two regular cospectral graphs and $H$ be any arbitrary regular graph, then

\begin{itemize}
    \item [(i)] $G\vee H$ and $G^{'}\vee H$ are $\mathcal{D}^S$-cospectral.
    \item [(ii)] $EDC(G)$ and $EDC(G^{'})$ are $\mathcal{D}^S$-cospectral.
\end{itemize}
\end{corollary}

\begin{corollary}
Let $G$ and $G^{'}$ be two $\mathcal{D}^S$-cospectral graphs, then 
\begin{itemize}
    \item [(i)] $G[K_2]$ and $G^{'}[K_2]$ are $\mathcal{D}^S$-cospectral.
    \item [(ii)]$D_2G$ and $D_2G^{'}$ are $\mathcal{D}^S$-cospectral.
\end{itemize}
\end{corollary}

A graph $G$ having all its $\mathcal{D}^S$-eigenvalues as integers is said to be \textit{distance Seidel integral} ($\mathcal{D}^S$-integral) graph. Clearly the graphs $K_{n}$ and $K_{n,n}$ are $\mathcal{D}^S$-integral. Following corollary's give some family of integral graphs formed by using different graph operations.
\begin{corollary}
Let $G= K_{n_1}$ and $G^{'}= K_{n_2}$, then $G\vee G^{'}$ is distance Seidel integral.
\end{corollary}
\begin{corollary}
For a distance Seidel integral graph $G$, its double graph $D_2G$ is also distance Seidel integral.
\end{corollary}
\begin{corollary}
Let $G$ be a distance integral graph, then $G\times K_2$ is distance Seidel integral.
\end{corollary}
\begin{corollary}
For a distance Seidel integral graph $G$, then $G[K_2]$ is distance Seidel integral.
\end{corollary}
\begin{corollary}
For a $k$-regular integral graph $G$, its extended double cover graph $EDC(G)$ is also distance Seidel integral.

\end{corollary}
\section{Conclusion}
In this article, a new graph matrix related to the Seidel matrix and distance matrix of the connected graph is defined. The relation between $\mathcal{D}^S$-eigenvalues with $\mathcal{D}$-eigenvalues and $A$-eigenvalues are established.  We characterized all the connected graphs with $\partial_{1}^{S}(G)= 3.$ Also, the bounds for the distance Seidel spectral radius and distance Seidel energy are obtained. The distance Seidel energy of different classes of graphs are calculated. In addition, the distance Seidel energy change of the complete bipartite graph due to the deletion of an edge is studied. Moreover, $\mathcal{D}^S$-spectra of graphs formed by using various graph operations are determined. In addition, different families of distance Seidel cospectral and distance Seidel integral graphs are obtained.
\section{Acknowledgements} 
The authors would like to thank the DST, Government of India, for providing support to carry out this work under the scheme 'FIST'(No.SR/FST/MS-I/2019/40). 
\section{Declarations}
 On behalf of all authors, the corresponding author states that there is no conflict of interest.
\bibliography{bib}
\bibliographystyle{plain}



\end{document}